\documentclass[10 pt, a4paper]{article}
\usepackage[T1]{fontenc}
\usepackage[utf8]{inputenc}
\usepackage[english]{babel}
\usepackage{mathtools}
\usepackage{amssymb}
\usepackage{amsmath}
\usepackage{amsthm}
\usepackage{geometry}
\usepackage{midpage}
\usepackage{graphicx}
\usepackage{tabularx}
\usepackage{array}
\usepackage{comment}
\usepackage{mathrsfs}
\usepackage{tikz}
\usepackage{tikz-cd}
\usepackage[bottom]{footmisc}
\usepackage{latexsym}
\usepackage{emptypage}
\usepackage{fancyhdr}
\usepackage{regexpatch}
\usepackage[hidelinks]{hyperref}

\title{Density of Noether-Lefschetz loci for surfaces in Fano and Calabi-Yau threefolds}
\author{Edoardo Mason}
\date{}

\begin{document}
\maketitle

\begin{abstract}
In this paper we provide applications of general results of Baldi-Klingler-Ullmo and Khelifa-Urbanik on the geometry of the Hodge locus associated to an integral polarized variation of Hodge structures to the case of Noether-Lefschetz loci for families of surfaces. In particular, we consider the family of surfaces in the linear system of a sufficiently ample line bundle on a smooth projective threefold $Y$, in case $Y$ is a Fano or a Calabi-Yau threefold, and we discuss the different behaviour of the union of the general, respectively exceptional, components of its Noether-Lefschetz locus.
\end{abstract}

\section{Introduction}

\subsection{Context}
It is an important problem in Hodge theory to understand the geometry of the Hodge locus associated with an integral polarized variation of Hodge structures ($\mathbb{Z}$-VHS from now on) on a smooth connected quasi-projective variety $S$, in particular it is interesting to have criteria to understand the density of families of components of the Hodge locus. Such a criterion is proved in Khelifa-Urbanik \cite{khur}.

This result fits into a more general line of work which has been developed in the last years, see for instance \cite{kli},\cite{klinglerhodge},\cite{bakker2020tame},\cite{bku1}. The essential point of view that has been adopted to work in this direction consists in seeing Hodge loci as intersection loci: indeed the period map $\Phi:S^{\text{an}} \to \Gamma \backslash D$ sends irreducible components of the Hodge locus to irreducible components of the intersection, inside $\Gamma \backslash D$, of the image of the period map with \emph{special} subvarieties of $\Gamma \backslash D$ which arise as quotients of period sub-domains of $D$ corresponding to Mumford-Tate groups that are smaller than the generic one. This perspective suggests to relate the geometry of special subvarieties of $S$ with the geometry of the corresponding intersection in $\Gamma \backslash D$, leading to the dichotomy between typical and atypical Hodge locus, discussed in Baldi-Klingler-Ullmo \cite{bku1}, and to the following conjecture on the behaviour of the two different parts of the Hodge locus.

\theoremstyle{plain}
\newtheorem{c11}{Conjecture}[section]
\begin{c11} \label{conj}
\emph{(\cite{bku1}, Conjecture 2.5 and 2.7)} Let $(\mathbb{V},\mathscr{F}^{\bullet})$ be an integral polarized variation of Hodge structures on a smooth connected quasi-projective variety $S$. Then:
\begin{itemize}
\item[1)] the typical Hodge locus for $(\mathbb{V},\mathscr{F}^{\bullet})$ is either empty or dense in $S^{\text{\emph{an}}}$ for the complex analytic topology, 
\item[2)] the atypical Hodge locus for $(\mathbb{V},\mathscr{F}^{\bullet})$ is algebraic, namely the set of atypical special subvarieties of $S$ has finitely many maximal elements for the inclusion.
\end{itemize}
\end{c11}

Results in the direction of this conjecture have been proven in the case of special subvarieties of positive period dimension, see in particular \cite{bku1} (Theorem 3.1) and \cite{khur}  (Theorem 1.9 and 3.5).

\subsection{The density criterion}
Let $S$ be a smooth connected quasi-projective variety and let $(\mathbb{V},\mathscr{F}^{\bullet})$ be a $\mathbb{Z}$-VHS on $S^{\text{an}}$, with generic Hodge datum $(\mathbf{G},D)$. For a strict Hodge sub-datum $(\mathbf{M},D_M)$ of $(\mathbf{G},D)$, we define the Hodge locus of type $\mathbf{M}$ as
\[
\text{HL}(S,\mathbb{V}^{\otimes},\mathbf{M}) = \{s \in S^{\text{an}}:\mathbf{MT}(\mathbb{V}_s) \subseteq g\mathbf{M}g^{-1} \text{ for some } g \in \mathbf{G}(\mathbb{Q})^{+} \} \subseteq \text{HL}(S,\mathbb{V}^{\otimes}),
\]
where for $s \in S^{\text{an}}$, $\mathbf{MT}(\mathbb{V}_s)$ is the Mumford-Tate group of the Hodge structure induced on the stalk $\mathbb{V}_s$.
The \emph{typical} Hodge locus of type $\mathbf{M}$ is defined to be the union $\text{HL}(S,\mathbb{V}^{\otimes},\mathbf{M})_{\text{typ}}$ of those components of $\text{HL}(S,\mathbb{V}^{\otimes},\mathbf{M})$ that are typical with respect to their generic Hodge datum. Then, Khelifa-Urbanik \cite{khur} proved the following:

\theoremstyle{plain}
\newtheorem{t11}[c11]{Theorem}
\begin{t11} \label{dense}
\emph{(\cite{khur}, Theorem 1.9)} Let $(\mathbb{V},\mathscr{F}^{\bullet})$ be a polarized $\mathbb{Z}$-VHS on a smooth connected quasi-projective variety $S$ with generic Hodge datum $(\mathbf{G},D)$. Furthermore, assume that the algebraic monodromy group of $\mathbb{V}$ is $\mathbf{H} = \mathbf{G}^{\text{\emph{der}}}$ and is $\mathbb{Q}$-simple. Let $(\mathbf{M},D_M)$ be a strict Hodge sub-datum of $(\mathbf{G},D)$ such that
\begin{equation} \label{eq:31}
\dim \Phi(S^{\text{\emph{an}}}) + \dim D_M - \dim D \geq 0.
\end{equation}
Then the Hodge locus of type $\mathbf{M}$ is analytically dense in $S^{\text{\emph{an}}}$.

Furthermore, if the inequality is strict, then the typical Hodge locus of type $\mathbf{M}$ is analytically dense.
\end{t11}

\subsection{Noether-Lefschetz loci}
This density criterion can be applied to study concrete geometric problems. In Khelifa-Urbanik \cite{khur}, it is used to give a new proof of the classical result of Ciliberto-Harris-Miranda \cite{chm} on the density of the (components of maximal codimension of the) Noether-Lefschetz locus for the universal family of smooth surfaces of degree $d \geq 5$ in $\mathbb{P}^3$. On the other hand, a known case of Conjecture \ref{conj} (\cite{bku1}, Theorem 3.1) allows to deduce results on the non density of families of components that have codimension strictly less than the maximum, as discussed in \cite{bku2}.

In this paper we further extend these applications, by applying the results of Khelifa-Urbanik (\cite{khur}, Theorem 1.9) and Baldi-Klingler-Ullmo (\cite{bku1}, Theorem 3.1) to other geometric contexts. First, for complete intersection surfaces in $\mathbb{P}^4$ we obtain:

\theoremstyle{plain}
\newtheorem{p11}[c11]{Proposition}
\begin{p11}
The typical part of the Noether-Lefschetz locus of the universal family of smooth complete intersection surfaces in $\mathbb{P}^4$ of bidegree $(d_1,d_2)$ with $2 \leq d_1 \leq 5$ and $d_2 \geq 3$ is analytically dense. On the other hand, the atypical part is algebraic.
\end{p11}

Moreover, we discuss the question of the density of the Noether-Lefschetz locus for the family of smooth surfaces in the linear system of a very ample line bundle $L$ on a smooth projective threefold $Y$, namely we consider the space $U_{Y,L}$ parametrizing smooth surfaces $X \subseteq Y$ such that $\mathscr{O}_Y(X) = L$ and the locus 
\[
\text{NL}(Y,L) = \{X \in U_{Y,L}:\text{Pic}(Y) \to \text{Pic}(X) \text{ is not surjective}\}.
\]
We will underline the relationship between the distinction between typical and atypical components of $\text{NL}(Y,L)$ and that between general and exceptional components, which only considers their codimension, and discuss the behaviour of these parts of $\text{NL}(Y,L)$ in the analyzed cases, thus illustrating in some geometric situations the Zilber-Pink paradigm introduced by Baldi-Klingler-Ullmo. Our aim here is to present some explicit results, in particular in the cases where $Y$ is a Fano or a Calabi-Yau threefold, for an abstract discussion see \cite{bku2} (Theorem 20).

\theoremstyle{plain}
\newtheorem{p12}[c11]{Proposition}
\begin{p12}
Let $Y$ be a smooth Fano threefold and $L$ be a very ample line bundle on $Y$. Then, for $d$ sufficiently large:
\begin{itemize}
\item[1)] the union of the typical components of $\text{\emph{NL}}(Y,L^d)$ is analytically dense in $U_{Y,L^d}$,
\item[2)] the union of the atypical components of $\text{\emph{NL}}(Y,L^d)$ is algebraic in $U_{Y,L^d}$.
\end{itemize}
\end{p12}

\theoremstyle{plain}
\newtheorem{p13}[c11]{Proposition}
\begin{p13}
Let $Y$ be a Calabi-Yau threefold with $\dim \text{\emph{Aut}}(Y) = 0$ and $L$ be a very ample line bundle on $Y$. Then, for $d$ sufficiently large:
\begin{itemize}
\item[1)] the union of the general components of $\text{\emph{NL}}(Y,L^d)$ is analytically dense in $U_{Y,L^d}$,
\item[2)] the union of the exceptional components of $\text{\emph{NL}}(Y,L^d)$ is algebraic in $U_{Y,L^d}$.
\end{itemize}
\end{p13}

We remark that the density of the Noether-Lefschetz locus for surfaces in a Calabi-Yau threefold was originally proved by Voisin \cite{voisindensite} (see also \cite{voisinintegral}) using the infinitesimal density criterion (\cite{voisinhodge}, Proposition 5.20).

\subsection{Notations}

Throughout this work, all algebraic varieties are defined over $\mathbb{C}$.

If $\mathbf{G}$ is a rational algebraic group, we denote by $\mathbf{G}^{\text{ad}}$ its adjoint group, by $\mathbf{G}^{\text{der}}$ its derived subgroup and by $\mathbf{G}(\mathbb{R})^{+}$ the identity connected component (for the real analytic topology) of $\mathbf{G}(\mathbb{R})$. Moreover, we define $\mathbf{G}(\mathbb{Q})^{+} = \mathbf{G}(\mathbb{R})^{+} \cap \mathbf{G}(\mathbb{Q})$.

\subsection{Acknowledgements}

Part of this work is included in the author's Master Thesis at Duisburg-Essen Universität. I want to express my gratitude to my advisor Dr. Carolina Tamborini for her constant support throughout this work and for all the interesting discussions on Hodge theory.

The author is currently supported by the project grant n. VR2023-03837.

\section{Background material}
We will review here some notions from the theory of variations of Hodge structures.

Recall that a Hodge structure on a rational vector space $V$ can be described as a real algebraic representation $\rho:\mathbf{S} \to \mathbf{GL}(V_{\mathbb{R}})$ of the Deligne torus. Let us fix such a representation $\rho$ and let $\mathbf{MT}(\rho)$ be its associated Mumford-Tate group, namely the smallest $\mathbb{Q}$-subgroup $\mathbf{M}$ of $\mathbf{GL}(V)$ such that $\rho$ factors through $\mathbf{M}_{\mathbb{R}}$. The associated Mumford-Tate domain is the orbit of $\rho$ in $\text{Hom}(\mathbf{S},\mathbf{MT}(\rho)_{\mathbb{R}})$ under the identity connected component $\mathbf{MT}(\rho)(\mathbb{R})^{+}$ of the real Lie group of $\mathbb{R}$-valued points of $\mathbf{MT}(\rho)$. We can summarize this situation in the following definition.

\theoremstyle{definition}
\newtheorem{d21}{Definition}[section]
\begin{d21}
\begin{itemize}
	\item[1)] A \emph{Hodge datum} is a pair $(\mathbf{G},D)$ where $\mathbf{G}$ is the Mumford-Tate group of some Hodge structure and $D$ is its associated Mumford-Tate domain.
	\item[2)] A morphism of Hodge data $(\mathbf{G},D) \to (\mathbf{G}',D')$ is a morphism of rational algebraic groups $\mathbf{G} \to \mathbf{G}'$ sending $D$ to $D'$.
	\item[3)] A Hodge sub-datum of $(\mathbf{G},D)$ is a Hodge datum $(\mathbf{G}',D')$ such that $\mathbf{G}'$ is a rational algebraic subgroup of $\mathbf{G}$ and the inclusion $\mathbf{G}' \to \mathbf{G}$ induces a morphism of Hodge data $(\mathbf{G}',D') \to (\mathbf{G},D)$. 
	\item[4)] A \emph{Hodge variety} is a quotient variety of the form $\Gamma\backslash D$ for some Hodge datum $(\mathbf{G},D)$ and some torsion-free arithmetic lattice $\Gamma \subseteq \mathbf{G}(\mathbb{Q})^{+}$, that is a torsion-free subgroup commensurable with $\mathbf{G}(\mathbb{Z})^{+}$.
\end{itemize}
\end{d21}

Now, let $S$ be a smooth connected quasi-projective variety over $\mathbb{C}$ and let $(\mathbb{V},\mathscr{F}^{\bullet})$ be an integral polarized variation of Hodge structures on $S^{\text{an}}$. It is a classical result due to Cattani-Deligne-Kaplan \cite{cdk} that the Mumford-Tate group of the induced Hodge structures on the stalks of $\mathbb{V}$ is locally constant outside a countable union of irreducible algebraic subvarieties of $S$: this union is the Hodge locus $\text{HL}(S,\mathbb{V}^{\otimes})$. Fixing a Hodge generic point $s \in S^{\text{an}} \smallsetminus \text{HL}(S,\mathbb{V}^{\otimes})$ and its associated (generic) Mumford-Tate group $\mathbf{G}$ we can construct the Hodge datum $(\mathbf{G},D)$, that we will call \emph{generic Hodge datum} of the given polarized $\mathbb{Z}$-VHS. 

If $Y$ is an irreducible algebraic subvariety of $S$, we can apply the same construction to the restriction of $\mathbb{V}$ to the smooth locus of $Y$ and obtain the generic Hodge datum $(\mathbf{G}_Y,D_{G_Y})$ of $Y$ for $(\mathbb{V},\mathscr{F}^{\bullet})$. Clearly, we can consider it as a Hodge sub-datum of $(\mathbf{G},D)$.

\theoremstyle{plain}
\newtheorem{p21}[d21]{Definition}
\begin{p21} \label{special}
An irreducible algebraic subvariety $Y$ of $S$ is special if and only if it is maximal for the inclusion among irreducible algebraic subvarieties of $S$ whose generic Mumford-Tate group is $\mathbf{G}_Y$.
\end{p21}

Clearly, $S$ itself is a special subvariety, while \emph{strict} special subvarieties are exactly the irreducible components of the Hodge locus.

Now recall that we can describe the given variation of Hodge structures by means of its period map $\Phi:S^{\text{an}} \to \Gamma\backslash D$, where $\Gamma$ is the monodromy group. We notice that we can characterise special subvarieties of $S$ as preimages under $\Phi$ of some specific analytic subvarieties of the Hodge variety $\Gamma\backslash D$, namely those which are Hodge varieties with respect to a Hodge sub-datum of $(\mathbf{G},D)$:

\theoremstyle{definition}
\newtheorem{d22}[d21]{Definition}
\begin{d22}
Given a Hodge variety $\Gamma\backslash D$, a \emph{special subvariety} of $\Gamma\backslash D$ is a subvariety of the form $\Gamma'\backslash D'$ for a Hodge sub-datum $(\mathbf{G}',D')$ of $(\mathbf{G},D)$ and $\Gamma' = \Gamma \cap \mathbf{G}'(\mathbb{Q})^{+}$.
\end{d22}

\theoremstyle{plain}
\newtheorem{p22}[d21]{Proposition}
\begin{p22} \label{special2}
\emph{(\cite{kli}, Proposition 3.20)} An irreducible algebraic subvariety $Y$ of $S$ is special if and only if $Y$ is an irreducible component of the preimage under $\Phi$ of a special subvariety of the Hodge variety $\Gamma\backslash D$.
\end{p22}

\theoremstyle{definition}
\newtheorem{d23}[d21]{Definition}
\begin{d23} \label{atypical}
A special subvariety $Y$ of $S$ with generic Hodge datum $(\mathbf{G}_Y,D_{G_Y})$ is said to be \emph{atypical} if $\Phi(S^{\text{an}})$ and $\Gamma_{\mathbf{G}_Y}\backslash D_{G_Y}$ intersect along $\Phi(Y^{\text{an}})$ with dimension bigger than expected, namely
\[
\text{codim}_{\Gamma\backslash D} \Phi(Y^{\text{an}}) < \text{codim}_{\Gamma\backslash D} \Phi(S^{\text{an}}) + \text{codim}_{\Gamma\backslash D} \Gamma_{\mathbf{G}_Y}\backslash D_{G_Y}.
\]
Otherwise, it is called \emph{typical}.
\end{d23}

\theoremstyle{remark}
\newtheorem{r21}[d21]{Remark}
\begin{r21}
This distinction can be found in Baldi-Klingler-Ullmo \cite{bku1} with the difference that special subvarieties of $S$ with image under $\Phi$ contained in the singular locus of $\Phi(S^{\text{an}})$ were considered always atypical. Following Khelifa-Urbanik \cite{khur} (Remark 2.8) we gave this definition since, being interested in density results for $\text{HL}(S,\mathbb{V}^{\otimes})$, we can always replace $S$ by a Zariski open subset, so we are free to assume that $\Phi(S^{\text{an}})$ is smooth.
\end{r21}

Then, the Hodge locus splits as
\[
\text{HL}(S,\mathbb{V}^{\otimes}) = \text{HL}(S,\mathbb{V}^{\otimes})_{\text{typ}} \cup \text{HL}(S,\mathbb{V}^{\otimes})_{\text{atyp}}
\]
where the typical Hodge locus $\text{HL}(S,\mathbb{V}^{\otimes})_{\text{typ}}$ (respectively $\text{HL}(S,\mathbb{V}^{\otimes})_{\text{atyp}}$) is the union of the \emph{strict} typical (respectively atypical) special subvarieties of $S$.

In the direction of Conjecture \ref{conj}, Baldi-Kligler-Ullmo proved the following:

\theoremstyle{plain}
\newtheorem{t21}[d21]{Theorem}
\begin{t21} \label{geomzp}
\emph{(\cite{bku1}, Theorem 3.1)} Let $(\mathbb{V},\mathscr{F}^{\bullet})$ be an integral polarized variation of Hodge structures on a smooth connected quasi-projective variety $S$. Let $Z$ be an irreducible component of the Zariski closure of the union of special atypical subvarieties of positive period dimension. Then:
\begin{itemize}
\item[1)] either $Z$ is a maximal atypical special subvariety,
\item[2)] or the generic adjoint Hodge datum $(\mathbf{G}_Z^{\text{\emph{ad}}},D_{G_Z})$ decomposes as a non trivial product $(\mathbf{G}',D') \times (\mathbf{G}'',D'')$, inducing (after replacing $S$ by a finite étale cover if necessary)
\[
\Phi|_{Z^{\text{\emph{an}}}} = (\Phi',\Phi''):Z^{\text{\emph{an}}} \to \Gamma_{\mathbf{G}_Z} \backslash D_{G_Z} = \Gamma' \backslash D' \times \Gamma'' \backslash D'' \subseteq \Gamma \backslash D,
\]
such that $Z$ contains a Zariski dense set of atypical special subvarieties for $\Phi''$ of zero period dimension.
\end{itemize}
\end{t21}

\theoremstyle{definition}
\newtheorem{d24}[d21]{Definition}
\begin{d24} \label{type}
Let $(\mathbf{M},D_M)$ be a strict Hodge sub-datum of $(\mathbf{G},D)$. The Hodge locus of type $\mathbf{M}$ is
\[
\text{HL}(S,\mathbb{V}^{\otimes},\mathbf{M}) = \{s \in S^{\text{an}}:\mathbf{MT}(\mathbb{V}_s) \subseteq g\mathbf{M}g^{-1} \text{ for some } g \in \mathbf{G}(\mathbb{Q})^{+} \} \subseteq \text{HL}(S,\mathbb{V}^{\otimes}).
\]
\end{d24}

Recall that the local system $\mathbb{V} \otimes_{\mathbb{Z}} \mathbb{Q}$ corresponds to a representation 
\[
\rho: \pi_1(S^{\text{an}},s_0) \to \mathbf{GL}(\mathbb{V}_{s_0} \otimes \mathbb{Q}),
\]
and that we define the algebraic monodromy group $\mathbf{H}$ of $\mathbb{V}$ as the connected component of the identity of the $\mathbb{Q}$-Zariski closure of $\rho(\pi_1(S^{\text{an}},s_0))$.

\theoremstyle{plain}
\newtheorem{t22}[d21]{Theorem}
\begin{t22} \label{dense}
\emph{(\cite{khur}, Theorem 1.9)} Let $(\mathbb{V},\mathscr{F}^{\bullet})$ be a polarized $\mathbb{Z}$-VHS on a smooth connected quasi-projective variety $S$. Furthermore, assume that the algebraic monodromy group of $\mathbb{V}$ is $\mathbf{H} = \mathbf{G}^{\text{\emph{der}}}$ and is $\mathbb{Q}$-simple. Let $(\mathbf{M},D_M)$ be a strict Hodge sub-datum of $(\mathbf{G},D)$ such that
\begin{equation} \label{eq:31}
\dim \Phi(S^{\text{\emph{an}}}) + \dim D_M - \dim D \geq 0.
\end{equation}
Then the Hodge locus of type $\mathbf{M}$ is analytically dense in $S^{\text{\emph{an}}}$.

Furthermore, if the inequality is strict, then the typical Hodge locus of type $\mathbf{M}$ is analytically dense.
\end{t22}

We refer to Khelifa-Urbanik \cite{khur} (Theorem 3.5) for a more general statement dealing with the case where the algebraic monodromy group is not $\mathbb{Q}$-simple.

\section{Density of Noether-Lefschetz loci}

In Khelifa-Urbanik \cite{khur} (section 1.2.2), Theorem \ref{dense} is used to deduce the density of the (components of maximal codimension of the) Noether-Lefschetz locus for the universal family of smooth surfaces of degree $d \geq 5$ in $\mathbb{P}^3$. 

Let us briefly review their strategy. It is a classical result of Lefschetz that a generic surface $X$ of degree $d \geq 5$ in $\mathbb{P}^3$ has Picard number $\rho(X) = 1$, namely the space of Hodge classes of type $(1,1)$ in $H^2(X,\mathbb{C})$ is one dimensional, generated by the first Chern class of $\mathscr{O}_X(1)$. If we let $U_d$ denote the space, open in $\mathbb{P}(H^0(\mathbb{P}^3,\mathscr{O}_{\mathbb{P}^3}(d)))$, parametrizing smooth degree $d$ surfaces in $\mathbb{P}^3$ and we let $f:\mathcal{X} \to U_d$ be its associated universal family, then we see that the Noether-Lefschetz locus
\[
\text{NL}(d) = \{s \in U_d:\rho(X_s) > 1\}
\]
is contained in the Hodge locus $\text{HL}(U_d,\mathbb{V}^{\otimes})$, for the integral polarized variation of Hodge structures supported on $\mathbb{V} = (R^2f_{*}\mathbb{Z})_{\text{prim}}$. Moreover, since the condition $\rho(X_s) > 1$ is equivalent to $H^2(X_s,\mathbb{Z})_{\text{prim}}$ containing a non zero Hodge class, we have $\text{HL}(U_d,\mathbb{V}^{\otimes}, \mathbf{M}) \subseteq \text{NL}(d)$, where $\mathbf{M}$ is the subgroup of the generic Mumford-Tate group that fixes a single Hodge vector. The analytic density of $\text{NL}(d)$ then follows from an application of Theorem \ref{dense} to the Hodge sub-datum $(\mathbf{M},D_M)$, since one can check that
\[
\dim \Phi(U_d) - \text{codim}_D D_M = \dim \Phi(U_d) - h^{2,0} > 0
\]
where $h^{2,0} = h^2(X,\mathscr{O}_X)$, for $X$ a smooth surface of degree $d$ in $\mathbb{P}^3$. More precisely, Theorem \ref{dense} gives the density of the typical components of the Noether-Lefschetz locus.

Let us clarify here the relationship between the dichotomy between typical and atypical special components of $\text{NL}(d)$ and the more classical one between \emph{general} and \emph{exceptional} components, as done in Baldi-Klingler-Ullmo \cite{bku2}. Notice first of all that for any irreducible component $Y$ of $\text{NL}(d)$ we have $\text{codim}_{U_d}Y \leq h^{2,0}$.

\theoremstyle{plain}
\newtheorem{r31}{Definition}[section]
\begin{r31} \label{defgen}
An irreducible component $Y$ of $\text{\emph{NL}}(d)$ is said to be general if $\text{\emph{codim}}_{U_d}Y = h^{2,0}$. Otherwise, it is said to be exceptional.
\end{r31}

\theoremstyle{plain}
\newtheorem{d31}[r31]{Proposition}
\begin{d31} \label{exc}
\emph{(\cite{bku2}, Proposition 17)}
\begin{itemize}
\item[1)] An irreducible component $Y$ of $\text{\emph{NL}}(d)$ is typical if and only if it is general and it has the expected generic Mumford-Tate group, namely isogenous to 
\[
\mathbf{G}_m \times \mathbf{SO}(2h^{2,0}, h^{1,1}_{\text{\emph{prim}}} - 1).
\]

\item[2)] If an irreducible component $Y$ of $\text{\emph{NL}}(d)$ is exceptional, then it is atypical.
\end{itemize}
\end{d31}

In particular, the previous argument implies the analytic density of the general components of $\text{NL}(d)$. On the other hand, an application of Theorem \ref{geomzp} (see \cite{bku2}, Theorem 15) implies that the union of atypical components of $\text{NL}(d)$ is algebraic, in particular this holds for the union of the exceptional components.

Let us consider now the case of smooth (scheme-theoretic) complete intersection surfaces in $\mathbb{P}^4$. First of all recall that the Noether-Lefschetz Theorem holds also in this case, except for the intersection of two quadrics (see Spandaw \cite{noether}). So, let us consider the space $U_{d_1,d_2}$ parametrizing smooth complete intersections of two hypersurfaces in $\mathbb{P}^4$ having degrees $d_1, d_2$ and define the Noether-Lefschetz locus
\[
\text{NL}(d_1,d_2) = \{s \in U_{d_1,d_2}:\rho(X_s) > 1\}.
\]

\theoremstyle{plain}
\newtheorem{p31}[r31]{Proposition}
\begin{p31} \label{complint}
The typical part of the Noether-Lefschetz locus of the universal family of smooth complete intersection surfaces in $\mathbb{P}^4$ of bidegree $(d_1,d_2)$ with $2 \leq d_1 \leq 5$ and $d_2 \geq 3$ is analytically dense. On the other hand, the atypical part is algebraic.
\end{p31}

\begin{proof}
First of all notice that, if we denote by $Q$ the polarization on $\mathbb{V} = (R^2f_{*}\mathbb{Z})_{\text{prim}}$ induced by cup-product, then a result of Beauville (\cite{beauville}, Theorem 5) ensures that the algebraic monodromy group $\mathbf{H}$ of $\mathbb{V}$ is the full orthogonal group $\text{Aut}(\mathbb{V}_s,Q_s)$, where we fixed some $s \in U_d$. Since we always have that $\mathbf{H}$ is a normal subgroup of the derived subgroup of the generic Mumford-Tate group $\mathbf{G}$ of $\mathbb{V}$ (see \cite{andre}, Theorem 1), this implies that $\mathbf{H} = \mathbf{G}^{\text{der}}$ and that it is $\mathbb{Q}$-simple. Thus we are in a situation where we can apply Theorem \ref{dense}.

Furthermore, by the local Torelli Theorem for complete intersections (see Flenner \cite{flenner}, Theorem 3.1) we have that the associated period map $\Phi$ satisfies
\[
\dim \Phi(U_{d_1,d_2}) = \dim U_{d_1,d_2} - \dim \text{PGL}_5(\mathbb{C}).
\]
Assuming without loss of generalities $d_1 \leq d_2$, one computes 
\[
\dim U_{d_1,d_2} = \left(\binom{d_1 + 4}{4} - 1\right) + \left(\binom{d_2 + 4}{4} - \binom{d_2 + 4 - d_1}{4} - 1\right).
\]
Now, if $(h^{2,0},h^{1,1},h^{0,2})$ denotes the Hodge numbers of $\mathbb{V}$ and $(\mathbf{M},D_M)$ is the Hodge sub-datum corresponding to the fixator of a single Hodge vector in $\mathbb{V}$, then $D_M$ is cut out in the full Mumford-Tate domain $D$ by $h^{2,0}$ equations, so $\dim D - \dim D_M = h^{2,0}$. 
To compute $h^{2,0}(X)$ for a complete intersection surface $X \subseteq \mathbb{P}^4$ given by the polynomials $f,g$, observe that one has a short exact sequence
\[
\begin{tikzcd}
0 \arrow[r] & \mathscr{O}_{\mathbb{P}^4}(-d_1 - d_2) \arrow[r] & \mathscr{O}_{\mathbb{P}^4}(-d_1) \oplus \mathscr{O}_{\mathbb{P}^4}(-d_2) \arrow[r] & \mathscr{I}_X \arrow[r] & 0
\end{tikzcd}
\]
where $\mathscr{I}_X \subseteq \mathscr{O}_{\mathbb{P}^4}$ is the ideal sheaf of $X$, the first map sends a local section $\sigma$ of $\mathscr{O}_{\mathbb{P}^4}(-d_1 - d_2)$ to $(\sigma g,- \sigma f)$ and the second map sends a local section $(\varphi,\psi)$ of $\mathscr{O}_{\mathbb{P}^4}(-d_1) \oplus \mathscr{O}_{\mathbb{P}^4}(-d_2)$ to $\varphi f + \psi g$. This yields
\[
0 \to H^3(\mathbb{P}^4,\mathscr{I}_X) \to H^4(\mathbb{P}^4,\mathscr{O}_{\mathbb{P}^4}(-d_1 - d_2)) \to H^4(\mathbb{P}^4, \mathscr{O}_{\mathbb{P}^4}(-d_1) \oplus \mathscr{O}_{\mathbb{P}^4}(-d_2)) \to H^4(\mathbb{P}^4,\mathscr{I}_X) \to 0.
\]
Now, considering the long exact sequence associated with the ideal sheaf sequence, we obtain $H^2(X,\mathscr{O}_X) \cong H^3(\mathbb{P}^4,\mathscr{I}_X)$ and $H^4(\mathbb{P}^4,\mathscr{I}_X) \cong H^3(X,\mathscr{O}_X) = 0$. Putting everything together we get
\begin{align*}
h^{2,0}(X) & = \dim H^2(X,\mathscr{O}_X) \\
& = \dim H^3(\mathbb{P}^4,\mathscr{I}_X) \\
& = \dim H^4(\mathbb{P}^4,\mathscr{O}_{\mathbb{P}^4}(-d_1 - d_2)) - \dim H^4(\mathbb{P}^4,\mathscr{O}_{\mathbb{P}^4}(-d_1)) - \dim H^4(\mathbb{P}^4,\mathscr{O}_{\mathbb{P}^4}(-d_2)) \\
& = \binom{d_1 + d_2 - 1}{4} - \binom{d_1 - 1}{4} - \binom{d_2 - 1}{4}.
\end{align*}
Thus, condition \ref{eq:31} for $(\mathbf{M},D_M)$ reads
\[
\binom{d_1 + 4}{4} + \binom{d_2 + 4}{4} - \binom{d_2 + 4 - d_1}{4} - 26 - \binom{d_1 + d_2 - 1}{4} + \binom{d_1 - 1}{4} + \binom{d_2 - 1}{4} > 0
\]
which is satisfied under the assumptions of the Proposition. We can therefore conclude the density of the typical part of $\text{NL}(d_1,d_2)$ by Theorem \ref{dense}.

Finally, the previous inequality implies that all components of $\text{NL}(d_1,d_2)$ have positive period dimension, hence Theorem \ref{geomzp}, together with the fact that the algebraic monodromy group of $\mathbb{V}$ is $\mathbb{Q}$-simple, implies that the union of atypical components of $\text{NL}(d_1,d_2)$ is algebraic, since option 2 of Theorem \ref{geomzp} would contradict the fact that $\mathbf{H}$ is $\mathbb{Q}$-simple.
\end{proof}

Let us now consider the following more general setting. Let $Y$ be a smooth projective threefold, $L$ be a very ample line bundle on $Y$ and $U_{Y,L}$ be the parameter space of smooth surfaces $X \subseteq Y$ such that $\mathscr{O}_Y(X) = L$. Moishezon \cite{moishezon} proved that a sufficient condition so that
\[
\text{Pic}(Y) \to \text{Pic}(X)
\]
is surjective for a general $X$ in $U_{Y,L}$ is that $h^{2,0}(X) > h^{2,0}(Y)$.
Notice that if $d$ is sufficiently large, then this condition is always satisfied by the line bundle $L^d$. Indeed, passing to the long exact cohomology sequence associated with
\[
\begin{tikzcd}
0 \arrow[r] & L^{-d} \arrow[r] & \mathscr{O}_Y \arrow[r] & \mathscr{O}_X \arrow[r] & 0,
\end{tikzcd}
\]
we obtain
\[\begin{tikzcd}
H^2(Y,\mathscr{O}_Y) \arrow[r] & H^2(X,\mathscr{O}_X) \arrow[r] & H^3(Y,L^{-d}) \arrow[r] & H^3(Y,\mathscr{O}_Y)
\end{tikzcd}
\]
and $h^3(Y,L^{-d}) = h^0(Y,L^d \otimes \omega_Y) \to \infty$ as $d \to \infty$, where $\omega_Y$ is the canonical bundle of $Y$.
Let us define the Noether-Lefschetz locus for the couples $(Y,L^d)$ as
\[
\text{NL}(Y,L^d) = \{X \in U_{Y,L^d}:\text{Pic}(Y) \to \text{Pic}(X) \text{ is not surjective}\}.
\]
We interpret $\text{NL}(Y,L^d)$ as a subset of the Hodge locus of the family of surfaces $X$ in $U_{Y,L^d}$ as done in the beginning of this section for surfaces in $\mathbb{P}^3$, in particular all irreducible components $Y$ of $\text{NL}(Y,L^d)$ have 
\[
\text{codim}_{U_{Y,L^d}}Y \leq h^{2,0}(X),
\]
thus we can give the same definition of general and exceptional components as in Definition \ref{defgen} and Proposition \ref{exc} still holds.

\theoremstyle{plain}
\newtheorem{p32}[r31]{Proposition}
\begin{p32} \label{fano}
Let $Y$ be a smooth Fano threefold and $L$ be a very ample line bundle on $Y$. Then, for $d$ sufficiently large:
\begin{itemize}
\item[1)] the union of the typical components of $\text{\emph{NL}}(Y,L^d)$ is analytically dense in $U_{Y,L^d}$,
\item[2)] the union of the atypical components of $\text{\emph{NL}}(Y,L^d)$ is algebraic in $U_{Y,L^d}$.
\end{itemize}
\end{p32}

\begin{proof}
We keep the same notation as in the proof of Proposition \ref{complint}, since we follow the same strategy. First of all notice that since we are interested in results for $d$ sufficiently large we can assume:
\begin{itemize}
\item[1)] The algebraic monodromy group of $\mathbb{V}$ coincides with $\mathbf{G}^{\text{der}}$ and is $\mathbb{Q}$-simple. Indeed by Carlson-Toledo \cite{monodromy} (Theorem 9.1), the image of the representation of $\pi_1(U_{Y,L^d})$ associated with $\mathbb{V}$ is Zariski dense in the group of automorphisms of $H^2(X,\mathbb{Q})_{\text{prim}}$ that respect the cup-product.
\item[2)] The period map $\Phi:U_{Y,L^d} \to \Gamma \backslash D$ has injective differential, in particular $\dim \Phi(U_{Y,L^d}) = \dim U_{Y,L^d} - \dim \text{Aut}(Y)$. This follows from Green \cite{green} (Theorem 0.1).
\end{itemize}
Since $Y$ is Fano we have $h^i(Y,\mathscr{O}_Y) = 0$ for $i \geq 1$ by the Kodaira vanishing Theorem, hence we can deduce from the ideal sheaf sequence for a general smooth surface $X$ in the linear system of $L^d$ that
\[
h^{2,0}(X) = h^3(Y,L^{-d}) = h^0(Y,L^d \otimes \omega_Y).
\]
By Iskovskih (\cite{iskovskih1978fano}, Proposition 1.3) we have $h^0(Y,\omega_Y^{-1}) \geq 4$, so we can write $\omega_Y = \mathscr{O}_Y(-D)$ for an effective divisor $D$. But then,  tensoring with $L^d$ the short exact sequence
\[
\begin{tikzcd}
0 \arrow[r] & \omega_Y \arrow[r] & \mathscr{O}_Y \arrow[r] & \mathscr{O}_D \arrow[r] & 0,
\end{tikzcd}
\]
we get the exact fragment
\[
\begin{tikzcd}
0 \arrow[r] & H^0(Y,L^d \otimes \omega_Y) \arrow[r] & H^0(Y,L^d) \arrow[r] & H^0(D,L^d \otimes \mathscr{O}_D) \arrow[r] & H^1(Y,L^d \otimes \omega_Y).
\end{tikzcd}
\]
Taking $d$ sufficiently large so that $H^1(Y,L^d \otimes \omega_Y) = 0$ we thus have
\[
h^0(Y,L^d) - h^0(Y,L^d \otimes \omega_Y) = h^0(D,L^d \otimes \mathscr{O}_D)
\]
which goes to infinity as $d \to \infty$, since the restriction of $L^d$ to $D$ is ample. Therefore, if $(\mathbf{M},D_M)$ is the Hodge sub-datum corresponding to the fixator of a single Hodge vector in $\mathbb{V}$, we have
\begin{align*}
\dim \Phi(U_{Y,L^d}) + \dim D_M - \dim D & = \dim U_{Y,L^d} -\dim \text{Aut}(Y) - h^{2,0}(X) \\
& = h^0(Y,L^d) - 1 - \dim \text{Aut}(Y) -  h^0(Y,L^d \otimes \omega_Y) > 0
\end{align*}
for $d$ sufficiently large. Hence Theorem \ref{dense} give the analytic density of the typical Hodge locus of type $\mathbf{M}$, in particular the union of the typical components of $\text{NL}(Y,L^d)$ is analytically dense in $U_{Y,L^d}$. Moreover, since the previous inequality implies that all components of $\text{NL}(Y,L^d)$ have positive period dimension, we can deduce from Theorem \ref{geomzp} that the union of atypical components of $\text{NL}(Y,L^d)$ is algebraic in $U_{Y,L^d}$.
\end{proof}

\theoremstyle{plain}
\newtheorem{p33}[r31]{Proposition}
\begin{p33} \label{cy}
Let $Y$ be a Calabi-Yau threefold with $\dim \text{\emph{Aut}}(Y) = 0$ and $L$ be a very ample line bundle on $Y$. Then, for $d$ sufficiently large:
\begin{itemize}
\item[1)] the union of the general components of $\text{\emph{NL}}(Y,L^d)$ is analytically dense in $U_{Y,L^d}$,
\item[2)] the union of the exceptional components of $\text{\emph{NL}}(Y,L^d)$ is algebraic in $U_{Y,L^d}$.
\end{itemize}
\end{p33}

\begin{proof}
The same considerations in the beginning of the proof of Proposition \ref{fano} hold in this case. To compute $h^{2,0}(X)$ for a general smooth surface in the linear system of $L^d$ we use the ideal sheaf sequence for $X$ and the well known Hodge numbers of a Calabi-Yau threefold $Y$ and we obtain
\[
h^{2,0}(X) = h^3(Y,L^{-d}) - h^3(Y,\mathscr{O}_Y) = h^0(Y,L^d) - 1.
\]
Therefore, keeping the same notations as in the previous proof,
\[
\dim \Phi(U_{Y,L^d}) + \dim D_M - \dim D = \dim U_{Y,L^d} - h^{2,0} = 0,
\]
hence Theorem \ref{dense} implies the analytic density of the Noether-Lefschetz locus for $(Y,L^d)$. On the other hand, since exceptional components of $\text{NL}(Y,L^d)$ have positive period dimension and are atypical, Theorem \ref{geomzp} implies that their union is algebraic. These two facts together imply that the union of the general components of $\text{NL}(Y,L^d)$ is analytically dense.
\end{proof}

\theoremstyle{remark}
\newtheorem{r32}[r31]{Remark}
\begin{r32}
Notice the subtle difference with the result of Proposition \ref{fano}. In the situation of Proposition \ref{cy} there could be atypical components of $\text{NL}(Y,L^d)$ which are general, hence have period dimension zero. It is still not known whether their union is not Zariski dense, since Theorem \ref{geomzp} only applies to components of positive period dimension. For the same reason we cannot establish whether the typical part of $\text{NL}(Y,L^d)$ is dense. 
\end{r32}

\theoremstyle{remark}
\newtheorem{r33}[r31]{Remark}
\begin{r33}
We remark that the density of the Noether-Lefschetz locus for surfaces in Calabi-Yau threefolds (proved originally by Voisin \cite{voisindensite}) has various geometric applications. It is used to show the nontriviality of the Abel-Jacobi map for Calabi-Yau threefolds (see \cite{voisindensite} and \cite{MR1741781}) and to prove the integral Hodge conjecture on Calabi-Yau threefolds (see \cite{voisinintegral}). Also the integral Hodge conjecture for curve classes on Fano fourfolds relies on this result (see \cite{MR2918165}).
\end{r33}

\section{Other applications}

\subsection{A family of hypersurfaces in $\mathbb{P}^5$}

Let us consider now the universal family $f:\mathcal{X} \to U_{n,d}$ of smooth hypersurfaces of degree $d$ in $\mathbb{P}^n$. We are interested in the Hodge locus of the polarized $\mathbb{Z}$-VHS associated with the primitive cohomology in middle degree, namely $\mathbb{V} = (R^{n-1}f_{*}\mathbb{Z})_{\text{prim}}$. Baldi-Klingler-Ullmo (\cite{bku1}, Corollary 1.6) proved that for $n \geq 4, d \geq 5$, $(n,d) \neq (5,5)$, the union of the components of $\text{HL}(U_{n,d},\mathbb{V}^{\otimes})$ of positive period dimension is algebraic: this follows from the fact that for such $n,d$ the level of this polarized $\mathbb{Z}$-VHS on $\mathbb{V}$ is $\geq 3$. Moreover, it is proved (\cite{bku1}, Theorem 3.3) that the typical Hodge locus for a polarized $\mathbb{Z}$-VHS of level $\geq 3$ is empty, hence it is conjectured that the whole Hodge locus $\text{HL}(U_{n,d},\mathbb{V}^{\otimes})$ is algebraic. Let us show here that in the case $n = d = 5$, excluded from this result, the typical Hodge locus is actually analytically dense. More precisely we prove:

\theoremstyle{plain}
\newtheorem{p41}{Proposition}[section]
\begin{p41}
The typical part of the \emph{generalized} Noether-Lefschetz locus 
\[
\text{\emph{NL}}(5,5) = \{s \in U_{5,5}: \text{\emph{rk}}(H^4(X_s,\mathbb{Z}) \cap H^{2,2}(X_s)) > 1\}
\]
of the universal family of smooth hypersurfaces of degree $5$ in $\mathbb{P}^5$ is analytically dense in $U_{5,5}$.
\end{p41}

\begin{proof}
First of all notice that for a general hypersurface $X_s$, $s \in U_{5,5}$, we indeed have $H^4(X_s,\mathbb{Z}) \cap H^{2,2}(X_s) \cong \mathbb{Z}$. This follows for instance from a result of Shioda (\cite{shioda}, Theorem 2.1) on the space of algebraic cohomology classes in  $H^4(X_s,\mathbb{Z}) \cap H^{2,2}(X_s)$ together with the fact that the Hodge conjecture is known for quintic fourfolds (see Da Silva Jr \cite{dasilva}, Corollary 2.20).

We can thus argue in the same spirit as the previous proofs. By Donagi \cite{donagi}, the period map is generically injective modulo the action of $\text{PGL}_6(\mathbb{C})$, thus 
\[
\dim \Phi(U_{5,5}) = \dim U_{5,5} - \dim \text{PGL}_6(\mathbb{C}).
\]
By Beauville (\cite{beauville}, Theorem 2) the algebraic monodromy group $\mathbf{H}$ of $\mathbb{V}$ coincides with the derived subgroup of the generic Mumford-Tate group and is $\mathbb{Q}$-simple. Consider the Hodge sub-datum $(\mathbf{M},D_M)$, where $\mathbf{M}$ is the fixator of a single Hodge vector, then 
\[
\dim D - \dim D_M \leq h^{0,4} + h^{1,3}.
\]
We easily compute $h^{4,0}(X) = 0$ for any $X$ in $U_{5,5}$ using the ideal sheaf sequence. To compute $h^{1,3}(X)$ notice that 
\[
h^{1,3}(X) = - 1 - \chi(X,\Omega_X^1)
\]
since for any hypersurface $X$ in $\mathbb{P}^n$ one has $H^q(X,\Omega_X^p) \cong H^q(\mathbb{P}^n,\Omega_{\mathbb{P}^n}^p) \cong \delta_{pq}\mathbb{C}$ for $p + q < n-1$. Furthermore, using the additivity of the Euler-Poincaré characteristic with the respect to the pull-back along the immersion $\iota:X \to \mathbb{P}^5$ of the Euler sequence for $\mathbb{P}^5$ and the conormal sequence
\[
\begin{tikzcd}
0 \arrow[r] & \mathscr{O}_{X}(-5) \arrow[r] &  \iota^{*}\Omega_{\mathbb{P}^5}^1 \arrow[r] & \Omega_X^1 \arrow[r] & 0
\end{tikzcd}
\]
we obtain 
\[
\chi(X,\Omega_X^1) = \chi(X,\mathscr{O}_X(-1)^{6}) - \chi(X,\mathscr{O}_X) - \chi(X,\mathscr{O}_X(-5)).
\]
Thus we compute
\[
\dim \Phi(U_{5,5}) + \dim D_M - \dim D \geq \dim U_{5,5} - \dim \text{PGL}_6(\mathbb{C}) - h^{1,3}(X) = 96
\]
and we can conclude by Theorem \ref{dense}.
\end{proof}

\subsection{Curves with non-simple Jacobian}

In this section we show how Theorem \ref{dense} easily implies a result of Colombo-Pirola (\cite{colpir}, Theorem 3) on the density of curves with non-simple Jacobian. Let $g \geq 4$ and consider the moduli space $\mathcal{M}_g$ of smooth projective curves of genus $g$. It is well known that the Zariski closure of the image of the Torelli map $j:\mathcal{M}_g \to \mathcal{A}_g$ to the moduli space $\mathcal{A}_g$ of principally polarized abelian varieties of dimension $g$ is Hodge generic (see Moonen-Oort \cite{moonenoort}, Remark 4.5) for the natural polarized $\mathbb{Z}$-VHS $\mathbb{V}$ induced by $j$. Notice that the algebraic monodromy group of $\mathbb{V}$ is $\mathbf{Sp}_{2g}$, in particular it is $\mathbb{Q}$-simple and coincides with the derived subgroup of the generic Mumford-Tate group.

Now let $1 \leq k \leq g-1$ be an integer and consider
\[
\mathcal{D}_{k} = \{C \in \mathcal{M}_g:J(C) \text{ contains an abelian subvariety of dimension $k$}\} \subseteq \mathcal{M}_g.
\]

\theoremstyle{plain}
\newtheorem{p42}[p41]{Proposition}
\begin{p42}
For $g \geq 4$ and $1 \leq k \leq 3$, $\mathcal{D}_k$ is analytically dense in $\mathcal{M}_g$.
\end{p42}

\begin{proof}
First of all notice that the condition of an abelian variety $A$ containing a proper abelian subvariety $B$ of dimension $k$ is equivalent to $A$ being isogenous to the product of $B$ with another subvariety of $A$, as follows from Poincaré reducibility Theorem. In other words, $\mathcal{D}_k$ is in correspondence with the set of intersections of $j(\mathcal{M}_g)$ with rational translates of the special subvariety $\mathcal{A}_k \times \mathcal{A}_{g-k}$ of $\mathcal{A}_g$. The condition
\[
\dim \mathcal{M}_g + \dim (\mathcal{A}_k \times \mathcal{A}_{g-k}) - \dim \mathcal{A}_g \geq 0
\]
of Theorem \ref{dense}
reads
\[
3g -3 + \frac{k(k+1)}{2} + \frac{(g-k)(g-k+1)}{2} - \frac{g(g+1)}{2} \geq 0,
\]
which simplifies to
\[
(3-k)g + k^2 - 3 \geq 0,
\]
which is clearly satisfied under the assumptions of the Proposition. We can thus conclude by Theorem \ref{dense}.
\end{proof}

\bibliographystyle{plain}
{\footnotesize
\bibliography{biblio_paper.bib}}
\end{document}